\documentclass[leqno, 11pt]{amsart}

\usepackage{lineno}
\usepackage{fullpage}
\usepackage{setspace}
\usepackage{tikz}
\usetikzlibrary{er}
\usetikzlibrary{arrows.meta}
\usepackage{amsmath}
\usepackage{amsfonts}
\usepackage{amssymb, mathrsfs}
\usepackage{verbatim,enumitem,graphics}
\usepackage{float}
\usepackage{microtype}
\usepackage{xcolor}
\usepackage[backref=page]{hyperref}
\hypersetup{%
    colorlinks,
    linkcolor={red!50!black},
    citecolor={green!50!black},
    urlcolor={blue!50!black}
}

\usepackage{dsfont}
\usepackage{array}

\usepackage[utf8]{inputenc}

\usepackage[abbrev,msc-links,backrefs]{amsrefs} 
\usepackage{doi}

\renewcommand{\PrintDOI}[1]{\doi{#1}}

\usepackage{fnpct}

\usepackage{MnSymbol}

\setenumerate{leftmargin=*} 

\usepackage{enumitem}
\usepackage[algo2e, english,algoruled,vlined,resetcount,linesnumbered]{
algorithm2e}
\usetikzlibrary{shapes,arrows,positioning}

\newtheorem{theorem}{Theorem}

\newtheorem{claim}[theorem]{Claim}

\newtheorem{construction}[theorem]{Construction}

\newtheorem{defn}[theorem]{Definition}

\newtheorem{fact}[theorem]{Fact}
\newtheorem{lemma}[theorem]{Lemma}

\newtheorem{problem}[theorem]{Problem}

\newtheorem{obs}[theorem]{Observation}
\numberwithin{equation}{section}
\numberwithin{theorem}{section}
\numberwithin{case}{section}

\def\F{\mathcal{F}}

\def \a{\alpha}
\def \e{\varepsilon}
\def \r{\gamma}

\def \bfi{\mathbf{i}}
\def \bfu{\mathbf{u}}
\def \bfv{\mathbf{v}}
\def \calP{\mathcal{P}}

\def \DPM{\textbf{DPM}}

\newcommand{\mb}[1]{\mathbb{#1}}

\newcommand{\nib}[1]{\noindent {\bf #1}}

\newcommand{\sub}{\subseteq}

\newcommand{\sm}{\setminus}

\newcommand{\es}{\emptyset}

\begin{document}


\title{Finding matchings in dense hypergraphs}
\author{Jie Han} 
\address{School of Mathematics and Statistics and Center for Applied Mathematics, Beijing Institute of Technology, Beijing, China, 100081}
\email{han.jie@bit.edu.cn}

\author{Peter Keevash}
\address{Mathematical Institute, University of Oxford, 
Woodstock Road, Oxford, OX2 6GG, UK.} 
\email{keevash@maths.ox.ac.uk}
\thanks{Supported by ERC Advanced Grant 883810.}
\date{}

\maketitle

\begin{abstract}
We consider the algorithmic decision problem
that takes as input an  $n$-vertex $k$-uniform 
hypergraph $H$ with minimum codegree at least $m-c$
and decides whether it has a matching of size $m$.
We show that this decision problem
is fixed parameter tractable with respect to $c$.
Furthermore, our algorithm not only decides the problem,
but actually either finds a matching of size $m$
or a certificate that no such matching exists.
In particular, when $m=n/k$ and $c=O(\log n)$, 
this gives a polynomial-time algorithm, that given
any $n$-vertex $k$-uniform hypergraph $H$ 
with minimum codegree at least $n/k-c$, 
finds either a perfect matching in $H$ 
or a certificate that no perfect matching exists.
\end{abstract}

\section{Introduction}

Matchings are fundamental objects in Graph Theory and have broad applications in other branches of Science and a variety of practical problems
(e.g.\ the assignment of graduating medical students to their first hospital appointments\footnote{In 2012, the Nobel Memorial Prize in Economics was awarded to Shapley and Roth ``for the theory of stable allocations and the practice of market design."}).
Applications of matchings in hypergraphs 
include the `Santa Claus' allocation problem~\cite{santa_claus};
they also offer a universal framework for many 
important combinatorial problems, e.g.\ the Existence Conjecture 
for designs (see \cite{Keevash_design, GKLO}) and
Ryser's conjecture \cite{ryser} on transversals in Latin squares.

This paper is concerned with the algorithmic question 
of finding a matching that is \emph{perfect},
meaning that it covers all vertices of the graph or hypergraph.
The graph case of this question is well understood:
Tutte's Theorem \cite{Tu47} gives necessary and sufficient conditions for a graph to contain a perfect matching, and Edmonds' Algorithm \cite{Edmonds} finds such a matching in polynomial time.
However, for hypergraphs it is a different story:
in fact, determining whether a $3$-uniform hypergraph contains a perfect matching was one of Karp's celebrated 21 NP-complete problems \cite{Karp}.
As the general problem is intractable (assuming P $\neq$ NP), 
it is natural to seek conditions that guarantee a perfect matching,
or at least make the existence question tractable.

\subsection{Perfect matchings under minimum degree conditions}
We start with some definitions that will be used throughout the paper.
Given $k\ge 2$, a $k$-uniform hypergraph (in short, \emph{$k$-graph}) $H=(V, E)$ consists of a vertex set $V$ and an edge set $E\subseteq \binom{V}{k}$, where every edge is a $k$-element subset of $V$. A \emph{matching} in $H$ is a collection of vertex-disjoint edges of $H$. A \emph{perfect matching} $M$ in $H$ is a matching that covers all vertices of $H$. We always assume that $k$ divides $n:=|V(H)|$, which is clearly a necessary condition for the existence of a perfect matching in $H$. 
{The perfect matching problem is also known as the \emph{exact cover} problem and can be viewed as a stricter form of the \emph{set cover} problem. Finding a matching of certain size is known as the \emph{$k$-set packing} problem: the decision problem for the existence of a given number of disjoint subsets from a family $\mathcal F$ of sets of size at most $k$. See also the survey~\cite{Hoffman2025}.}

{We will consider the algorithmic effect 
of minimum degree conditions defined as follows. 
For $S \sub V(H)$ the \emph{neighbourhood} of $S$ is 
$N_H(S) := \{ T \sub V(H) \sm S: S \cup T \in E(H) \}$,
and the \emph{degree} of $S$ is $\deg_{H} (S) = |N_H(S)|$;
the subscript $H$ is omitted if it is clear from the context. 
The \emph{minimum $d$-degree $\delta _{d} (H)$} 
of $H$ is the minimum of $\deg_{H} (S)$ over all $d$-vertex sets $S$ in $H$.}

{We refer to $\delta_{k-1} (H)$ as the \emph{minimum codegree} of $H$.
This is the usual parameter considered for minimum degree results 
on hypergraph matchings. It is also natural to consider other $d$,
although the case $d=0$, which concerns the size of $H$, 
leads to a rather trivial question, as it is not hard to see
that the largest size is attained by a complete hypergraph 
on $n-1$ vertices together with an isolated vertex.}

R\"odl, Ruci\'nski and Szemer\'edi \cite{RRS09} determined 
the sharp minimum codegree condition to ensure a perfect matching 
in an $n$-vertex $k$-graph for large $n$ and all $k\ge 3$.
The quantity is $n/2-k+C$, where $C\in\{3/2, 2, 5/2, 3\}$ depends on the values of $n$ and $k$, and the extremal examples are as follows.

\begin{construction}
\label{con13}
Given $n\ge k\ge 3$, let $\mathcal{H}_{n,k}$ be the collection of $k$-graphs $H$ such that there is a partition of $V(H)= X\cup Y$ such that $n/k-|X|$ is odd and all edges of $H$ intersect $X$ in an odd number of vertices. 
\end{construction}

Any matching $M$ in $H\in \mathcal H_{n,k}$ 
cannot be perfect; indeed, if it were we would have
$\sum_{e\in M} (|e\cap X|-1) = |X|-n/k$,
which is impossible, as the left hand side is even
but the right hand side is odd.
There is a large literature on minimum degree 
conditions that force a perfect matching, see
\cite{AFHRRS, CzKa, FrKu_mat, HPS, Khan2, Khan1, KO06mat, KOT, MaRu, Pik, RRS06mat, RRS09, TrZh12, TrZh13, TrZh15} and the surveys \cite{RR, zsurvey},
yet there are still many open problems, such as determining 
the minimum 1-degree condition that forces a perfect matching.



\subsection{Algorithms}

Let $\DPM_k(n, m)$ be the decision problem of determining 
whether an $n$-vertex $k$-graph $H$ with $\delta_{k-1}(H)\ge m$ 
contains a perfect matching. When can $\DPM_k(n, m)$ be decided
in polynonmial time? 

$\DPM_k(n,0)$ is well-known to be NP-complete for constant $k\ge 3$ due to Exact 3-cover (X3C) being one of the famous Karp’s 21 NP-complete problems.
The result of \cite{RRS09} mentioned above
shows that the decision problem is trivial for $m \ge n/2 - k + 3$
(then there is a perfect matching iff $k \mid n$).
Szyma\'nska \cite{Szy13} proved that if $\delta<1/k$ then
$\DPM_k(n,0)$ is polynomial-time reducible to $\DPM_k(n,\delta n)$,
and so $\DPM_k(n,\delta n)$ is NP-complete. Karpi\'nski, Ruci\'nski and Szyma\'nska \cite{KRS10} showed that there exists $\e>0$ such that $\DPM_k(n,(1/2-\e)n)$ is in P 
and posed the question of determining 
for the complexity of $\DPM_k(n,\delta n)$ 
for $\delta\in [1/k, 1/2)$. 

This was resolved for $\delta>1/k$
by Keevash, Knox and Mycroft \cite{KKM_abs, KKM13},
who showed that the problem is in P 
and also gave a polynomial-time algorithm
that either finds a perfect matching or a certificate 
that none exists. Their proof was very long and technical, 
and left open the `threshold case' $\delta=1/k$,
which poses additional challenges 
(discussed in Section \ref{sec:over}).
Han \cite{Han14_poly} completely resolved the question 
of \cite{KRS10} by showing that $\DPM_k(n,\delta n)$ 
is in P for $\delta\in [1/k, 1/2)$. 
His proof was much simpler than that in \cite{KKM13}
(it relied on some theory from~\cite{KKM13}
and also developed a lattice-based absorbing method 
which has found many other applications)
but his result only concerned the decision problem,
and left open the constructive problem,
i.e.\ finding a perfect matching or a certificate 
that none exists.
{For a summary of these results, see Table~\ref{table}. As in this paper we consider $k$, the uniformty of the hypergraph, as a constant, we use $O_k(\cdot)$ notation when the implicit multiplicative constant is a function of $k$; on the other hand, we also use $O(\cdot)$ when the implicit multiplicative constant is an absolute constant.
}
\begin{table}[h!]
\centering
\caption{Known results on $\DPM_k(n,m)$. The
last column indicates whether the corresponding algorithm is constructive and outputs a perfect matching.}
\label{table}
\begin{tabular}{|c|c|c|c|} 
\hline
\textbf{Reference} & \textbf{Result} & \textbf{Running time} &\textbf{Finds PM} \\ 
\hline
\cite{Karp} & $\DPM_k(n,0)\in $ NP-C & N/A & N/A  \\ 
\hline
\cite{Szy13} & $\forall$ fixed $\varepsilon> 0$, $\DPM_k(n,n/k-n^{\e})\in $ NP-C & N/A & N/A  \\ 
\hline
\cite{KRS10} & $\exists$ fixed $\varepsilon>0$ s.t. $\DPM_k(n,(1/2-\e)n)\in $ P & Implicit & yes  \\ 
\hline
\cite{KKM13} & $\forall$ fixed $\varepsilon>0$, $\DPM_k(n,(1/k+\e)n)\in $ P & $O_k(n^{3k^2-7k+1}) + O_\e(1)$ & yes \\ 
\hline
\cite{Han14_poly} & $\DPM_k(n,n/k)\in $ P & $O_k(n^{3k^2-5k})$ & no \\ 
\hline
\cite{HaKe20} & $\forall$ fixed $c\ge 0$, $\DPM_k(n,n/k-c)\in $ P & $O_k(n^{k\cdot\max\{c,4^k\}})$ & yes \\ 
\hline
This work & $\forall$ $c=O_k(\log n)$, $\DPM_k(n,n/k-c)\in $ P & $O_k(n^{2^{k+1}k^2}+ 2^{O_k(c)})$ & yes \\ 
\hline
\end{tabular}
\end{table}

\subsection*{Main results}
The first contribution of this paper 
is the following result.

\begin{theorem}\label{thm:main}
Let $k\ge 3$ and $H$ be an $n$-vertex $k$-graph 
with $\delta_{k-1}(H)\ge n/k-c$ for some $c>0$. 
Then there is an algorithm running in time
$O_k(n^{2^{k+1}k^2}+ 2^{O_k(c)})$
that finds a perfect matching in $H$ 
or a certificate that none exists.
In particular, if $c = O_k(\log n)$
then $\DPM_k(n,n/k-c)$ is in P.
\end{theorem}

For example, a certificate from Theorem~\ref{thm:main} can be a $k$-graph $H_0$ as described in Construction~\ref{con13} which itself does not admit a perfect matching (which can be checked in polynomial-time) and includes $H$ as a subgraph.
We include detailed discussions in Section 5.

In the conference version of this paper~\cite{HaKe20},
we proved Theorem~\ref{thm:main} with a worse running time, namely, $O_k(n^{k\max\{4^k, c\}})$,
thereby obtaining a polynomial-time algorithm only when $c$ is a constant.
Our new improvement lies in the proof of the extremal case (Lemma~\ref{lem:t-d}), 
where the brute force search used in~\cite{HaKe20} is replaced by
Theorem~\ref{thm:Fellows} below, which is proved by 
the colour-coding technique of Alon, Yuster and Zwick~\cite{AYZ95}.


{Our main result deals with matchings of arbitrary size, and therefore, is a special case of the \textsc{$k$-Set Packing} problem.
We first review the following theorem of Han~\cite{Han14_mat}.
\begin{theorem}
\cite{Han14_mat}
Given $k\in \mathbb N$, there exists $n_0\in \mathbb N$ such that the following holds for all $n\ge n_0$.
If $H$ is a $k$-graph with $\delta_{k-1}(H)\ge m$ and $m < n/k$, then $H$ contains a matching of size $m$. 
\end{theorem}
In a similar fashion as Theorem~\ref{thm:main}, one may investigate the following parameterized problem.
\medskip
\begin{center}
\fbox{
    \begin{minipage}{0.9\textwidth}
        \textsc{$k$-Graph Matching Above Minimum Codegree}
        \hfill \textbf{Parameter:} \( c \) \\
        \textbf{Input:} \( k \)-uniform \( n \)-vertex hypergraph \( H \) with $\delta_{k-1}(H)\ge m$, an integer \( c \geq 1 \). \\
        \textbf{Problem:} Find a matching of size \(\min\{m + c, \lfloor \frac{n}{k} \rfloor\}\) in \( H \) or report that it does not exist.  
    \end{minipage}
}
\end{center}
\medskip
A result of Fellows, Knauer, Nishimura, Ragde, Rosamond, Stege, Thilikos, and Whitesides~\cite{FellowsFPT}~solved the $m=0$ case.
\begin{theorem}
\cite{FellowsFPT}
\label{thm:Fellows}
The problem \textsc{$k$-Graph Matching Above Minimum Codegree} with $m=0$ admits an algorithm with running time $O(n^{k}+ 2^{O(kc)})$.
\end{theorem}
We address this problem by showing the following result.
}

\begin{theorem}\label{coro}
Let $c, k, n\ge 3$ and $0<m\le n/k$.
Let $H$ be an $n$-vertex $k$-graph 
with $\delta_{k-1}(H)\ge m-c$.
Then there is an algorithm running in time
$O_k(n^{2^{k+1}k^2+k}+ 2^{O_k(c)}n^k)$
that finds a matching of size $m$ in $H$ 
or a certificate that none exists.
That is, the problem \textsc{$k$-Graph Matching Above Minimum Codegree} admits an algorithm with running time $O_k(n^{2^{k+1}k^2+k}+ 2^{O_k(c)}n^k)$.
\end{theorem}

The study of such algorithmic versions of classical problems
in extremal graph theory was recently proposed by Fomin, Golovach, Sagunov and Simonov~\cite{Fomin22},
who obtained corresponding algorithmic versions of some extremal results
for cycles in graphs, such as Dirac's theorem.

We note that Theorem \ref{thm:main} is a special case ($m=n/k$) of Theorem~\ref{coro}.
In fact, we will now show that the general statement
is a simple consequence of this case.

\begin{proof}[Proof of Theorem~\ref{coro} assuming Theorem~\ref{thm:main}]
We start by considering the case of $t:= \frac{n-mk}{k-1}\in \mathbb Z$.
Note that $m+t= (n+t)/k$, so $n+t\in k\mathbb N$.
Consider any $n$-vertex $k$-graph $H$ with $\delta_{k-1}(H)\ge m-c$.
We form a new $k$-graph $H'$ with $|V(H')|=n+t$ 
by adding a set $A$ of $t$ (new) vertices to $V(H)$ 
and adding all $k$-vertex sets intersecting $A$ to $E(H)$.
We note that $\delta_{k-1}(H')\ge m-c+t$.
Furthermore, writing $\nu$ for the matching number,
{we have $\nu(H) \ge m$ if and only if $ \nu(H')\ge m + t=(n+t)/k$}, 
as any matching in $H$ can be greedily extended
to a matching in $H'$ with $t$ extra edges,
and from any matching in $H'$ we can obtain a matching in $H$
by deleting at most $t$ edges that intersect $A$.
In particular, the existence of or finding matchings
of size $m$ in $H$ is algorithmically equivalent
to the same problem for perfect matchings in $H'$.
This can be solved by the algorithm given by Theorem~\ref{thm:main}, 
which finds a perfect matching in $H'$ or a certificate that none exists.
The running time is
\[
O_k((n+t)^{2^{k+1}k^2}+ 2^{O_k(c)}) = O_k(n^{2^{k+1}k^2}+ 2^{O_k(c)}).
\]
This proves the result if $t:= \frac{n-mk}{k-1}\in \mathbb Z$.
For the general case, we choose $s \le k-1$ such that 
$n_0:=n-s$ has $t:= \frac{n_0-mk}{k-1}\in \mathbb Z$.
We apply the above method to each of the $\tbinom{n}{s}$ 
possible $k$-graphs $H_0$ on $n_0$ vertices
obtained by deleting $s$ vertices from $H$,
noting that $\delta_{k-1}(H_0)\ge \delta_{k-1}(H)-s \ge m-c-(k-1)$.
\end{proof}


\section{Preliminaries} \label{sec:over}

In this section we give an overview of our proof.
\subsection{Extremal versus non-extremal}
In the next section we describe the algorithm referred
to in Theorem \ref{thm:main} and reduce its proof 
of correctness to two theorems (proved in the two
subsequent sections) that respectively handle
the `non-extremal' and `extremal' cases for $H$.
To explain this distinction, we consider the following
construction that appears naturally
around the codegree threshold $n/k$.

\begin{construction}[Space Barrier]\label{con:sb}
Let $V$ be a set of size $n$ 
and fix $S\subseteq V$ with $|S|<n/k$.
Let $H$ be the $k$-graph on $V$ 
whose edges are all $k$-sets that intersect $S$.
\end{construction}

We note that the minimum codegree of $H$ is $|S|$ 
and any matching in $H$ has at most $|S|<n/k$ edges,
so cannot be perfect. On the other hand,
Han \cite{Han14_mat} showed that any $n$-vertex $k$-graph
with $\delta_{k-1}(H)\ge n/k-1$ contains a matching 
of size $n/k-1$, thus determining the tight codegree condition
for a matching that is just one edge short of being perfect.
This rather surprising phenomenon indicates that the key
issue for whether there is a perfect matching near the codegree
threshold $n/k$ is whether $H$ is close to a space barrier,
or equivalently, whether $H$ has an independent set 
of size about $n-n/k$: we say $H$ is \emph{$\e$-extremal} 
if it has an independent set of size $(1-\e)\frac{k-1}{k}n$.

Our strategy of separating the non-extremal and extremal cases
follows that of Han \cite{Han14_poly}, and indeed several aspects 
of his proof carry over to our setting (despite our weaker assumption
$\delta_{k-1}(H)\ge n/k-c$, which is significantly more challenging
to work with, as it does not rule out space barriers).
However, the crucial difficulty that prevented Han from finding
a perfect matching (as opposed to just testing for its existence)
is the algorithmic intractability of finding (or even approximating)
the largest independent set (this is NP-hard even in graphs, 
let alone hypergraphs).

The starting point for the new approach in this paper is 
to observe that this issue can be avoided via an algorithm of
Han~\cite{Han14_mat} that either finds an almost perfect matching 
or a large independent set. In the latter case, we can algorithmically
extend to a large \emph{maximal} independent set, which may not be
of maximum size, but nevertheless gives us enough power 
to analyse the extremal case with some additional arguments.
A final contribution of this paper (see the concluding remarks)
is an algorithmic reduction of the perfect matching problem  
(see Problem~\ref{prob}) of independent interest.

\subsection{Divisibility and lattices}
In this subsection we introduce some key definitions and tools in the proofs, which are motivated by the divisibility construction (Construction~\ref{con13}).

Recall every $k$-graph in the family $\mathcal H_{n,k}$ (defined in Construction~\ref{con13}) contains no perfect matching, which is certified by a biparition of the vertices together with the restriction on the edges.
This divisibility construction plays an essential role for determining the exact minimum codegree condition guaranteeing perfect matching.
Moreover, similar ``divisibility'' constructions exist for an arbitrary number of parts (see~\cite{KKM13} for more examples).
Motivated by the divisibility barriers, Keevash, Knox and Mycroft \cite{KKM13} introduced the following concepts to capture the divisibility property of $H$.

\begin{defn} \label{def:L}
\cite{KKM13}
Let $H=(V, E)$ be a $k$-graph and $\calP=\{V_1, V_2, \dots, V_d\}$ 
be a partition of $V$. The index vector of a set $S\subseteq V$
is $\bfi_{\calP}(S)=(|S\cap V_1|, \dots, |S\cap V_d|) \in \mb{Z}^d$.
Given $\mu>0$, let $I_{\calP}^\mu(H)$ denote the set of all 
$k$-vectors $\bfi\in \mathbb{Z}^d$ such that 
at least $\mu |V|^k$ edges $e\in H$ have $\bfi_{\calP}(e)=\bfi$,
and let $L_{\calP}^{\mu}(H)$ denote the lattice in $\mathbb{Z}^d$ (additive subgroup of $\mathbb Z^d$)
generated by $I_{\calP}^\mu(H)$.
\end{defn}

Note that for a lattice $L$ (an additive subgroup of $\mathbb Z^d$) if $\bfi_{\calP}(e)\in L$ for every $e\in E(H)$ but $\bfi_{\calP}(V(H))\notin L$, then $H$ contains no perfect matching.
Indeed, for any matching $M$ in $H$, we have $\bfi_{\calP}(V(M)) = \sum_{e\in M} \bfi_{\calP}(e)\in L$, implying that $V(M)\neq V(H)$.
This proves the following fact.
{
\begin{fact}
Let $L$ be the lattice in $\mathbb{Z}^d$ generated by members of $I_{\calP}(H)$.
If $H$ contains a perfect matching, then $\bfi_{\calP}(V(H))\in L$.
\end{fact}
}
{However, as noted in~\cite{KKM13}, the condition $\bfi_{\calP}(V(H))\in L$ is not sufficient, and it is intuitively rather weak, as for each $\bfi\in L$ there may only be one edge of $H$ with index vector $\bfi$, and if two such edges intersect then they cannot both be used in a matching.} This motivates the definition of $I_{\calP}^{\mu}(H)$, where a vector $\bfi$ is included only if $H$ contains many edges with index vector $\bfi$.
Now the idea is that $L_{\calP}^{\mu}(H)$ captures 
the robust `divisibility' constraints in $H$ and that $H$ should not have any divisibility obstruction 
to a perfect matching if it is possible to delete a small matching
so that the index vector of the uncovered set is in $L_{\calP}^{\mu}(H)$.

We will also use the reachability methods 
introduced by Lo and Markstr\"om \cite{LM2, LM1}.

\begin{defn} \label{def:R}
Let $H$ be an $n$-vertex $k$-graph.
We say that two vertices $u$ and $v$ are \emph{$(\beta, i)$-reachable} 
in $H$ if there are at least $\beta n^{i k-1}$ $(i k-1)$-sets $S$ such that both $H[S\cup \{u\}]$ and $H[S\cup \{v\}]$ have perfect matchings. 
We say that $U \sub V(H)$ is \emph{$(\beta, i)$-closed in $H$} 
if any two vertices $u,v\in U$ are $(\beta, i)$-reachable in $H$.
\end{defn}


\section{The algorithm}

In this section we state our algorithm and prove our main theorem
assuming two theorems (concerning the non-extremal and extremal cases)
whose proofs will be given in the following two sections.

The first ingredient of our algorithm is
Procedure ListPartitions from \cite[Section 2]{KKM13},
which can test\footnote{
See~\cite[Lemma 2.2]{KKM13} 
and note that the proof is valid assuming 
$\delta_{k-1}(H) \ge \omega n$
for any fixed $\omega$ and $n>n_0(k,\omega)$.} 
in time $O_k(n^{k+1})$ whether $H\in \mathcal H_{n,k}$
(defined in Construction~\ref{con13}).
If $H\in \mathcal H_{n,k}$ then the procedure finds
a partition certifying that $H$ does not have a perfect matching,
so we reduce to the case $H\notin \mathcal H_{n,k}$.

The following statement
is a simplified form of \cite[Lemma 2.5]{Han14_poly}.
Throughout this paper, $x\ll y$ means that for any $y> 0$ there 
exists $x_0> 0$ such that for any $x< x_0$ the subsequent statement holds.

\begin{lemma}
\cite[Lemma 2.5]{Han14_poly}
\label{lem:PL}
Let $1/n \ll \beta \ll \r\ll 1/k$, where $k \ge 3$ is an integer.
Then for each $k$-graph $H$ on $n$ vertices 
with $\delta_{k-1}(H)\ge n/k-\r n$, 
Procedure LatticePartition in time $O_k(n^{2^{k-1}k+1})$ outputs
a partition $\calP=\{V_1, \dots, V_d \}$ such that 
each $V_i$ is $(\beta, 2^{k-1})$-closed in $H$
and has $|V_i|\ge n/k-2\r n$.
\end{lemma}

Note that we will always have $d\le k$ in the lemma above due to the lower bound on the size of each part $V_i$.

The following theorem handles the non-extremal case:
given the partition $\calP$ from Lemma \ref{lem:PL},
it provides an algorithm that finds a perfect matching
or outputs an independent set witnessing 
that $H$ is in the extremal case.

\begin{theorem}\label{thm2}
Suppose $1/n \ll \beta, \mu \ll \r \ll 1/k$.
Let $H$ be a $k$-graph on $n$ vertices 
such that $\delta_{k-1}(H)\ge n/k - \r n$. 
Suppose $\calP$ is a partition of $V(H)$ found by Lemma \ref{lem:PL}.
Suppose there is a matching $M_1$ with $|M_1| \le k$
such that $\bfi_{\calP}(V(H)\setminus V(M_1))\in L_{\calP}^{\mu}(H)$.
Then there is an algorithm that finds in time $O_k(n^{2^{k+1}k^2})$
a perfect matching in $H$ 
or an independent set in $H$ of size $(1-5k\r)\frac{k-1}{k}n$. 
\end{theorem}

The extremal case is handled by the following theorem
(recall that $H$ is {$\e$-extremal} if $V(H)$ 
contains an independent subset 
of size at least $(1-\e)\frac{k-1}k n$).

\begin{theorem}\label{thm3}
Assume $1/n\ll \e \ll 1/k$ and $c\le \e n/k$. 
Let $H$ be a $k$-graph on $n$ vertices 
such that $\delta_{k-1}(H)\ge n/k - c$. 
Suppose $H\notin \mathcal H_{n,k}$ and $H$ is $\e$-extremal. 
Then there is an algorithm Procedure~\ref{alg:determine} that finds
in time $O_k(n^{4k-4}+ 2^{O(kc)})$
a perfect matching in $H$ 
or a certificate that none exists.
\end{theorem}

Now we are ready to state our main algorithm (see Procedure~\ref{alg:main}).

\begin{procedure}[h]
  \caption{PerfectMatching()}
  \label{alg:main}
  \SetKwInOut{Input}{Data}
  \SetKwInOut{Output}{Output}
  \Input{an $n$-vertex $k$-graph $H$ such that $\delta_{k-1}(H)\ge n/k-c$.}
  \Output{a perfect matching in $H$ or a certificate that none exists.} 
Choose constants $1/n_0\ll \beta, \mu \ll\r\ll 1/k$ \;
  \If{$n\le n_0$}
  {Examine every set of $n/k$ edges in $H$, and halt with appropriate output.}
  \If{$c\ge 5\gamma n$}{Solve the problem by the subroutine in Theorem~\ref{thm:Fellows}, and halt with appropriate output.}
  \If{Procedure ListPartitions finds that $H\in \mathcal H_{n,k}$ }{
  Output the certifying partition and halt. 
  }
  Apply {Procedure LatticePartition} (Lemma~\ref{lem:PL}) to find a partition $\calP$ of $V(H)$ and $L_{\calP}^{\mu}(H)$ \;
      \eIf {there is a matching $M_1$ of size at most $k$ with
      $\bfi_{\calP}(V(H)\setminus V(M_1))\in L_{\calP}^{\mu}(H)$} {
    Find matchings $M_0$ and $M_2$ as in the proof of Theorem \ref{thm2};
    
    \eIf{Lemma~\ref{thm:next} applied on $H-(V(M_0\cup M_1\cup M_2))$ outputs a matching $M_3$}{
    Use Theorem~\ref{thm2} to find a perfect matching given $M_0, M_1, M_2$ and $M_3$; Output this perfect matching and halt.
    }{
    Apply {Procedure~\ref{alg:determine}} with $\e=5k\gamma$ and halt with appropriate output.
    }}{Output $\calP$ as a certificate of ``no perfect matching'' and halt.} 
\end{procedure}

Now we give an overview of our algorithm, Procedure~\ref{alg:main}.
The algorithm first chooses constants $n_0\in \mathbb N$, $\beta, \mu, \gamma$ satisfying the relations on line 1.
{In particular, they are \textit{all} functions of $k$.}
Given an input $k$-graph $H$ with appropriate minimum codegree condition, {the algorithm first deals with two trivial cases: if $n\le n_0$, then a trivial brute-force algorithm solves the problem in constant time; if $c\ge 5\gamma n$, then the subroutine in Theorem~\ref{thm:Fellows} can find a perfect matching or conclude its non-existence in time $O(n^k+2^{O(n)})=2^{O(c/\gamma)}=2^{O_k(c)}$.}
Then the subroutine Procedure ListPartitions is used to determine whether $H$ belongs to $\mathcal{H}_{n,k}$, with a vertex partition certifying the non-existence of perfect matchings in $H$.
If $H\notin\mathcal{H}_{n,k}$, then Procedure LatticePartition is used to find a partition $\mathcal P$ that captures the reachability information of vertices in $H$.
On line 7, the algorithm exhaustively checks the existence of a matching $M_1$ of size at most $k$ with $\bfi_{\calP}(V(H)\setminus V(M_1))\in L_{\calP}^{\mu}(H)$, which indeed is another \textit{necessary} condition for the existence of perfect matching in $H$.
In particular, if such $M_1$ does not exist, then the partition $\calP$ together with the non-existence of $M_1$ forms a checkable certificate of the non-existence of perfect matchings.
If such $M_1$ exists, then the algorithm constructs two useful matchings $M_0$ and $M_2$, with the property that if $H-V(M_0\cup M_1\cup M_2)$ contains a matching that leaves out at most $k^2/\gamma$ vertices, then a perfect matching can be constructed and outputted.
The algorithm then applies the subroutine in Lemma~\ref{thm:next} to $H-V(M_0\cup M_1\cup M_2)$; if it finds a suitable matching $M_3$ then a perfect matching is constructed, otherwise the subroutine outputs an independent set of size $(1-2k\gamma)(k-1)n/k$, from where another key subroutine Procedure~\ref{alg:determine} takes over.
With the existence of such a large independent set, $H$ exhibits a strong structure that will be further analyzed by Procedure~\ref{alg:determine}, and several substructures (also certifying by certain vertex partitions) are tested for non-existence of perfect matchings.
If no such substructure is found, then the algorithm constructs a perfect matching.

\begin{figure}
\centering
\resizebox{0.9\textwidth}{!}{
	\begin{tikzpicture}[
		node distance=2.5cm,
		decision/.style={diamond, draw, aspect=2,text width=5em, align=flush center, inner sep=0.3mm},
			decision1/.style={diamond, draw, aspect=2, text width=14em, align=flush center, inner sep=0.3mm, minimum height=1.5 },
			process1/.style={rectangle, draw,  minimum height=4em, text width=14em, align=center},
			process2/.style={rectangle, draw,  minimum height=4em, text width=14.5em, align=center},
		process/.style={rectangle, draw, text width=15em, align=flush center},
		terminal1/.style={ellipse, draw, text width=8em, align=center},
		terminal/.style={ellipse, draw, text width=5em, align=center},
			terminal2/.style={ellipse, draw, text width=6em, align=center},
		arrow/.style={->, >=stealth', ultra thick}
		]
		
		\node (start) [terminal1] {Start:\\ 
			Input an 
			$n$-vertex 
			$k$-graph $H$ with $\delta_{k-1}(H)\ge n/k-c$};
		
		\node (n0) [decision, below of=start,yshift=-0.2cm] {$n \leq n_0$?};
		
		\node (brute) [terminal, left of=n0, xshift=-3.5cm] {
			Solve via brute-force algorithm
		};
		
		\node (c) [decision, below of=n0] {$c \geq 5\gamma n$?};
		
		\node (thm) [terminal, left of=c, xshift=-3.5cm] {Solve via Theorem 1.4
		};
		
		\node (list) [decision, below of=c, yshift=-0.2cm] {Does Listpartitions find an $H \in \mathcal{H}_{n,k}$?};
		
		\node (cert) [terminal, left of=list, xshift=-3.5cm] {Output the certifying partition
		};
		
		\node (lattice) [process, below of=list,yshift=-0.2cm] {Procedure LatticePartition   finds $\mathcal{P}$ and $L_\mathcal{P}^\mu(H)$ (by Lemma 3.1) };
		
		\node (match) [decision1, right of=lattice,xshift=6cm] {
Does there exist a matching $M_1$ in $H$ such that 
$|M_1|\le k$ and $i_\mathcal{P}(V(H)\setminus V(M_1))
\in L_\mathcal{P}^\mu(H)$?};

		\node (no) [terminal1, right of=match, xshift=5.2cm] {Output $\mathcal{P}$ as a certificate of “no perfect matching”
		};

		\node (absorb) [process2, above of=match,yshift=0.7cm] {Find  matchings $M_0$ and $M_2$ \\(as in the proof of  Theorem 3.2)};
		
		\node (lemma) [decision1, above of=absorb,yshift=0.7cm] {Does Lemma 4.1 applied on  $H - (V(M_0\cup M_1\cup M_2))$  output a matching $M_3$?};
		
		\node (perfect) [terminal1, right of=lemma, xshift=5.2cm] {Use Theorem 3.2 to find a  perfect matching given $M_0,M_1, M_2,$ and $M_3$
		};
		
		\node (ext) [terminal2, above of=lemma, yshift=1.4cm] {Apply PerfectMatchingEXT 
		};

		\draw [arrow] (start) -- (n0);
		\draw [arrow] (n0) -- node[above] {Yes} (brute);
		\draw [arrow] (n0) -- node[left] {No} (c);
		\draw [arrow] (c) -- node[above] {Yes} (thm);
		\draw [arrow] (c) -- node[left] {No} (list);
		\draw [arrow] (list) -- node[above] {Yes} (cert);
		\draw [arrow] (list) -- node[left] {No} (lattice);
		\draw [arrow] (lattice) -- (match);
		\draw [arrow] (match) -- node[right] {Yes} (absorb);
		\draw [arrow] (match) -- node[above] {No} (no);
		\draw [arrow] (absorb) -- (lemma);
		\draw [arrow] (lemma) -- node[above] {Yes} (perfect);
		\draw [arrow] (lemma) -- 
		node[right,align=left]{No: \\
			$\exists$ a large independent set
		} (ext);

	\end{tikzpicture}
}
    \caption{A diagram view for the main algorithm, Procedure \ref{alg:main}}
\end{figure}

We conclude this section by showing correctness of
this algorithm, thus proving Theorem~\ref{thm:main},
assuming Theorems \ref{thm2} and \ref{thm3},
which will be proved in the following two sections.
Suppose $L$ is an edge-lattice in $\mathbb{Z}^{|\calP|}$, where $\calP$ is a partition of a set $V$, then the \emph{coset group} of $(\calP, L)$ is $G=G(\calP, L)=L_{\max}^{|\calP|}/L$, where $L_{\max}^d=\{ x\in \mathbb{Z}^d: k \text{ divides }\sum_{i\in [d]} x_i \}$.

\begin{proof}[Proof of Theorem~\ref{thm:main}]
We show correctness of Procedure~\ref{alg:main}.
{As discussed above, 
we may assume that $n\ge n_0$ and $c<5\gamma n$, as otherwise a trivial subroutine or the subroutine in Theorem~\ref{thm:Fellows} solves the problem accordingly, in constant time or in time $2^{O_k(c)}$.}
Moreover, Procedure ListPartitions from \cite[Section 2]{KKM13},
tests in time $O_k(n^{k+1})$ whether $H\in \mathcal H_{n,k}$
and if so outputs a partition 
certifying that $H$ does not have a perfect matching.
Thus we can assume $H\notin \mathcal H_{n,k}$.

Next we apply Lemma \ref{lem:PL} 
which finds $\calP = (V_1,\dots,V_d)$ 
in time $O_k(n^{2^{k-1}k+1})$;
we note that $d \le k$.
Then we test each set of at most $k$ edges in $H$
(of which there are $O_k(n^{k^2})$) to see if it is 
a matching $M_1$ satisfying
$\bfi_{\calP}(V(H)\setminus V(M_1))\in L_{\calP}^{\mu}(H)$\footnote{To check if a given vector $\bfv$ is in a lattice $L$, one can take an integer basis of $L$ and use Gaussian Elimination to solve for the (unique) real coefficients of $\bfv$ in this basis. Then $\bfv\in L$ if and only if all the coefficients are integers. The running time is a function of $k$.}.
If we find any such $M_1$ then we use it to apply
Theorem \ref{thm2}, which finds a perfect matching
or an independent set that can be used to apply Theorem \ref{thm3},
which in turn finds a perfect matching
or a certificate that none exists.
We only need to apply Theorems \ref{thm2} and \ref{thm3}
for one such $M_1$ if it exists, so the running time
is
\[
O_k(n^{2^{k+1}k^2}+ O_k(n^{4k-4}+ 2^{O(kc)}) = O_k(n^{2^{k+1}k^2}+ 2^{O(kc)}).
\]

To complete the proof, it remains to show that if there is
a perfect maching $M$ then some such $M_1$ exists.
To see this we argue similarly to \cite{KKM13}.
First we note the following property of $I_{\calP}^{\mu}(H)$
that follows easily from the minimum codegree condition of $H$:
for any $\mathbf{v} \in \mb{Z}^d$ with non-negative coordinates
summing to $k-1$ there is some $i\in [d]$ 
such that $\mathbf{v}+\bfu_i\in I_{\calP}^{\mu}(H)$.
Indeed, consider all $(k-1)$-sets with index vector $\bfv$, by the minimum codegree condition $H$ contains at least $(n/k-2\gamma n-k)^{k-1}\cdot \delta_{k-1}(H)/k! \ge (1/k-3\gamma)^k n^k/k!$ edges whose index vector equals $\bfv+\bfu_i$ for $i\in [d]$; then by averaging, there exists $i\in [d]$ such that the number of edges with index vector $\bfv+\bfu_i$ is at least $\mu n^k$.
Then the proof of \cite[Lemma 6.4]{KKM13} shows that
the coset group $G$ of $L_{\calP}^{\mu}(H)$ in the lattice
$\{\mathbf{v} \in \mb{Z}^d: k \mid \sum v_i \}$
has size $|G| \le d \le k$.
Now we repeatedly apply the pigeonhole principle to reduce $M$ 
to $M_1$ as in \cite[Proposition 6.10]{KKM13}.
We start with $M_1=M$ and note that as $M$ is perfect we have
$\bfi_{\calP}(V(H)\sm V(M))=\mathbf{0}\in L_{\calP}^{\mu}(H)$. 
While $|M_1|>k$ we consider any edges $e_1,\dots,e_k$ in $M_1$
and the partial sums $\sum_{j=1}^{i} \bfi_{\calP}(e_j)$ 
for $0 \le i \le k$. By the pigeonhole principle, 
some two of these sums lie in the same coset 
of $L_{\calP}^{\mu}(H)$, namely, there exist $0\le i_1< i_2\le k$ such that
\[
\sum_{j=i_1+1}^{i_2} \bfi_{\calP}(e_j)= \sum_{j=1}^{i_2} \bfi_{\calP}(e_j) - \sum_{j=1}^{i_1} \bfi_{\calP}(e_j)\in L_{\calP}^{\mu}(H).
\] 
So we can delete $e_{i_1+1},\dots, e_{i_2}$ from $M_1$ while preserving
$\bfi_{\calP}(V(H)\sm V(M_1)) \in L_{\calP}^{\mu}(H)$. 
We terminate with $|M_1| \le k$, as required.
\end{proof}

\section{The non-extremal case}

In this section we prove Theorem \ref{thm2},
which finds a perfect matching or a large independent set,
thus establishing correctness of Procedure \ref{alg:main}
in the non-extremal case. We adapt (and simplify) 
the approach of Han \cite{Han14_poly} 
via lattice-based absorption
and also incorporate a derandomisation argument
of Garbe and Mycroft \cite{GaMy} 
so that we can find a perfect matching
(not just test for its existence).

\subsection{Almost perfect matching or large independent set}

The key idea of proof via absorption is that it simplifies
the problem of finding a perfect matching to that
of finding an almost perfect matching. Accordingly,
we start by showing how to find an almost perfect matching
or a large independent set. The following lemma
is essentially \cite[Lemma 1.6]{Han14_mat};
the proof is algorithmic, although this is not made explicit,
so for the convenience of the reader we do this here.

\begin{lemma} \label{thm:next}
Suppose that $1/n\ll \r \ll 1/k$ and $k \mid n$. 
Let $H$ be a $k$-graph on $n$ vertices 
with $\delta_{k-1}(H)\ge n/k - \r n$. 
Then in time $O_k(n^{k+1})$ we can find 
either a matching that leaves at most $k^2/\gamma$ vertices uncovered, 
or an independent set of size $(1-2k\gamma)\frac{k-1}{k}n$.
\end{lemma}

\begin{proof}
Consider any matching $M=\{e_1, e_2, \dots, e_m\}$ in $H$. 
Let $V'$ be the set of vertices covered by $M$ 
and let $U$ be the set of vertices which are not covered by $M$. 
Assume that $|U|> k^2/\r$ and $U$ is an independent set 
(otherwise we can trivially enlarge $M$).
We will show that we can find either a matching of size $m+1$, 
or an independent set of size $(1-2k\gamma)\frac{k-1}{k}n$.

We arbitrarily partition all but at most $k-2$ vertices of $U$ 
into disjoint $(k-1)$-sets $A_1, \dots, A_t$ 
where $t=\lfloor \frac{|U|}{k-1}\rfloor>\frac{k}{\r}$.
Let $D$ be the set of vertices $v\in V'$ such that 
$\{v\}\cup A_i\in E(H)$ for at least $k$ choices of $A_i$.
First we consider the case that there is some $i^* \in [m]$
with $|e_{i^*} \cap D| \ge 2$. We fix distinct $x,y$ in $e_{i^*} \cap D$
and apply the definition of $D$ to pick distinct $A_i, A_j$ 
such that $\{x\}\cup A_i$ and $\{y\}\cup A_j$ are edges.
Then we can enlarge $M$ by replacing $e_{i^*}$ 
by $\{x\}\cup A_i$ and $\{y\}\cup A_j$.
Thus we may assume each $|e_i\cap D|\le 1$.

Next we show that $|D|\ge (\frac 1k -2\r)n$. 
As $\delta_{k-1}(H)\ge n/k - \r n$ and $U$ is independent,
\[ t\left(\frac 1k -\r \right)n
\le \sum_{i=1}^t \deg(A_i)
\le |D| t +n\cdot k < t(|D|+\r n),
\]
so we have
\[|D| > \left(\frac 1k -2\r \right)n.\]

Let $V_D:=\bigcup\{e_i: e_i\cap D\neq \emptyset\}$. 
Note that $|V_D\sm D|=(k-1)|D|\ge (k-1)(\frac 1k -2\r)n$. 
If $V_D \sm D$ is independent then we can output
it to conclude the proof. Thus we can assume
there is some edge $e_0 \in H[V_D\sm D]$.
Suppose that $e_0$ intersects $e_{i_1}, \dots, e_{i_l}$ in $M$ 
for some $l\le k$ and write $\{v_{i_j}\}= e_{i_j}\cap D$ 
for all $j\in [l]$. By definition of $D$, 
we can greedily pick $A_{i_1}, \dots, A_{i_l}$ 
such that  $\{v_{i_j}\}\cup A_{i_j} \in E(H)$ for all $j\in [l]$.
Then we can enlarge $M$ by replacing
the edges $e_{i_1}, \dots, e_{i_l}$ 
by $e_0$ and $\{v_{i_j}\}\cup A_{i_j}$ for $j\in [l]$. 

For the running time, note that there are at most $n/k$ iterations
in which we enlarge the matching.
In each iteration we find $D$ by calculating 
the degree of each vertex $v$ in time $O_k(n^{k-1})$.
Then we check if $V_D\sm D$ is an independent set 
in time $O_k(n^k)$. The overall running time is $O_k(n^{k+1})$.
\end{proof}

\subsection{An absorbing lemma}

Now we implement the absorption part of the algorithm, Lemma~\ref{lem:abs},
by derandomising a lemma from \cite{Han14_poly}.
For the statement, we recall 
Definitions \ref{def:L} and \ref{def:R},
and say that a set $T$
\emph{absorbs} or \emph{is absorbing for} a set $S$
if both $H[T]$ and $H[T\cup S]$ contain perfect matchings.

\begin{lemma}[Absorbing Lemma]\label{lem:abs}
Let $1/n \ll 1/c' \ll \beta, \mu \ll 1/k, 1/t, 1/C'$.
Suppose that $\calP=\{V_1, \dots, V_d\}$ is a partition of $V(H)$ 
such that each $V_i$ is $(\beta, t)$-closed in $H$.
Then in time $O(n^{4tk^2})$ we can find a family $\F_{abs}$ 
of at most $c'\log n$ disjoint $tk^2$-sets 
such that $H[V(F)]$ contains a perfect matching 
for all $F\in \F_{abs}$ and every $k$-vertex set $S$ 
with $\bfi_{\calP}(S)\in I_{\calP}^{\mu} (H)$ 
has at least $C'$ absorbing $t k^2$-sets in $\F_{abs}$. 
\end{lemma}

The proof uses the following lemma 
which is \cite[Claim 3.6]{Han14_poly}.

\begin{lemma}\label{clm:abs}
Suppose $V_i$ is $(\beta, t)$-closed in $H$ for all $i\in [d]$.
Then any $k$-set $S$ with $\bfi_{\calP}(S)\in I_{\calP}^{\mu}(H)$ 
has at least $\frac{\mu \beta^k}{2^{k+1}} n^{t k^2}$ absorbing $t k^2$-sets.
\end{lemma}

For the derandomisation we use the following lemma of 
Garbe and Mycroft~\cite[Proposition 4.7, Procedure SelectSet]{GaMy} 
(see also Karpi\'nski, Ruci\'nski and Szyma\'nska~\cite{KRS10hc}).

\begin{lemma}\label{lem:deran}
Fix constants $\beta> \tau>0$ and integers $m, M, N$ and $r\le N$ such that $r$ and $N$ are sufficiently large, and that $M\le (1/8) \exp(\tau^2r/(3\beta))$.
Let $U$ and $W$ be disjoint sets of sizes $|U|=M$ and $|W|=N$.
Let $G$ be a graph with vertex set $U\cup W$ such that $G[U]$ is empty, $G[W]$ has precisely $m$ edges, and $\deg_G(u)\ge \beta N$ for every $u\in U$.
Then in time $O(N^4+MN^3)$ we can find an independent set $R\subseteq W$ in $G$ such that $(1-\nu)r\le |R|\le r$ and $|N_G(u)\cap R| \ge(\beta - \tau - \nu)r$ for all $u\in U$, where $\nu = 2m r/N^2$. 
\end{lemma}

We are now ready to prove the absorbing lemma.

\begin{proof}[Proof of Lemma \ref{lem:abs}]
We will apply Lemma~\ref{lem:deran} 
to the graph $G$ with parts 
$U = \{S: \bfi_{\calP}(S)\in I_{\calP}^{\mu} (H)\}$ 
and $W = \{T \sub V(H): |T|=tk^2\}$,
where $T,T'$ in $W$ are adjacent iff $T \cap T' \ne \es$,
and $S \in U$ and $T \in W$ are adjacent iff $T$ absorbs $S$.
In the notation of Lemma~\ref{lem:deran}
we have $N=\binom{n}{tk^2}$, $M\le \binom n k$
and $m=|E(G[W])|\le tk^2 \binom{n}{tk^2-1}\binom{n}{tk^2} 
= t^2k^4N^2/(n-tk^2+1)$.
We let $\beta':=\mu \beta^k/2^{k+1}$, 
$\tau := \beta'/3$ and $r=c'\log n$.
Then by Lemma~\ref{clm:abs}, 
$\deg_G(u)\ge \beta' N$ for every $u\in U$, and 
\[
\exp\left( \frac{\tau^2 r}{3\beta'} \right) = \exp \left( \frac{\beta' c'\log n}{27} \right) \ge 8\binom nk \ge 8M,
\]
as $n$ is large enough.
Thus, by Lemma~\ref{lem:deran}, in time $O(N^4+M N^3) = O(n^{4tk^2})$, 
we can find a set $R\subseteq W$ which is independent in $G$ 
with $(1-\tau)r\le |R|\le r$ 
and $|N_G(u)\cap R| \ge(\beta' - \tau - \nu)r$ for all $u\in U$, where
\[
\nu = \frac{2m r}{N^2} \le \frac{2 t^2k^4N^2 c'\log n}{(n-tk^2+1)N^2} = \frac{2 t^2k^4 c'\log n}{n-tk^2+1} < \beta'/3.
\] 
Note that $R$ consists of disjoint $tk^2$-sets of $H$
by definition of $G[W]$.
We now remove $tk^2$-sets in $R$ that do not have a perfect matching, 
and denote the resulting family of $tk^2$-sets by $\F_{abs}$.
Thus $|\F_{abs}|\le c'\log n$, 
each member of $\F_{abs}$ has a perfect matching, 
and every $k$-vertex set $S$ with 
$\bfi_{\calP}(S)\in I_{\calP}^{\mu} (H)$ 
has at least $\beta' r/3 \ge \beta' c'\log n/3 \ge C'$ 
absorbing $t k^2$-sets in $\F_{abs}$, as $n$ is large enough.
\end{proof}

\subsection{Proof of Theorem \ref{thm2}}

%

We start the proof by fixing a parameter hierarchy
$1/n \ll 1/c'\ll \beta, \mu \ll \r, 1/C \ll 1/k$,
such that $C$ has the following property: for any  
$I \sub \{ \bfv \in \mb{Z}^d: \sum_i |v_i| \le k \}$,
where $d \le k$, and $\bfu$ in the integer span of $I$
with $\sum_i |u_i| \le k^2(1+\r^{-1})$, we can write
\begin{equation}\label{eq:Cmax}
\bfu=\sum_{\bfv\in I}a_{\bfv}(\bfu)\bfv. 
\end{equation}
where $a_{\bfv}(\bfu)\in \mathbb{Z}$
and $|a_{\bfv}(\bfu)| \le C$ for each $\bfv\in I$.

Let $H$ be a $k$-graph on $n$ vertices 
such that $\delta_{k-1}(H)\ge n/k - \r n$.
Let $\calP=\{V_1,\dots, V_d\}$ be the partition 
found by Lemma \ref{lem:PL};
we note that $d\le k$ by Lemma \ref{lem:PL}.
Write
\[C':=(C+k+2)k+k/\gamma,
\quad I=I_{\calP}^{\mu}(H), 
\quad \text{ and } \ L=L_{\calP}^{\mu}(H).
\]
Let $M_1$ be a matching of size at most $k$ 
such that $\bfi_{\calP}(V(H)\setminus V(M_1))\in L$.
We first apply Lemma \ref{lem:abs} to $H$ with $t=2^{k-1}$,
obtaining a family $\F_{abs}$ of $2^{k-1}k^2$-sets 
with $|\F_{abs}|\le c'\log n$ such that every $S$ 
with $\bfi_{\calP}(S)\in I$ has at least $C'$ 
absorbing $2^{k-1} k^2$-sets in $\F_{abs}$.
We obtain $\F_{0}$ from $\F_{abs}$ 
by removing any sets that intersect $V(M_1)$,
of which there are at most $|V(M_1)|\le k^2$;
then every $S$ with $\bfi_{\calP}(S)\in I$ 
has at least $C' - k^2$ absorbing sets in $\F_{0}$.
We let $M_0$ be a perfect matching on $V(\F_0)$ 
consisting of perfect matchings on each member of $\F_0$. 

Next, we greedily find a matching $M_2$ 
in $V(H)\setminus V(M_{0}\cup M_1)$ which contains $C$ edges $e$ 
with $\bfi_{\calP}(e)=\bfi$ for every $\bfi\in I$. 
This  is possible because $H$ contains at least $\mu n^k$ edges 
for each $k$-vector $\bfi\in I$, whereas 
$|V(M_0\cup M_{1}\cup M_2)|
\le 2^{k-1} k^2 c'\log n + k^2 + k C \tbinom{k+d-1}{k} < \mu n$, 
as $n$ is large enough.

Let $H':=H[V(H)\setminus V(M_0\cup M_{1} \cup M_2)]$. 
Note that $|V(H')|\ge n - \mu n$ as $n$ is large, and
\[
\delta_{k-1}(H') \ge \delta_{k-1}(H) - \mu n \ge n/k - 2\r n \ge (1/k - 2\r)|V(H')|.
\]
We apply Lemma~\ref{thm:next} to $H'$ with parameter $2\r$ in place of $\r$.
If Lemma~\ref{thm:next} finds an independent set of size $(1-2k(2\gamma))\frac{k-1}{k}|V(H')|$, then we can output this independent set and halt, as 
\[
(1-4k\r)\frac{k-1}{k}|V(H')| \ge (1-4k\r)\frac{k-1}{k} (n-\mu n)\ge (1-5k\r)\frac{k-1}{k}n.
\]
Thus we can assume that Lemma~\ref{thm:next} finds a matching $M_3$ 
covering all but a set $S_0$ 
of at most $k^2/(2\r)\le k^2/\r$ vertices of $V(H')$.

Recall that $\bfi_{\calP}(V(H)\setminus V(M_1))\in L$. 
By the definition of $M_2$, 
we have $\bfi_{\calP}(V(H)\setminus V(M_1\cup M_2))\in L$, 
which implies 
$\bfi_{\calP}(S_0) + \sum_{e \in M_{0}\cup M_3} \bfi_{\calP}(e) \in L$. 
As in the previous section,
by the pigeonhole principle we can find in time $O_k(n)$ edges 
$e_1, \dots, e_{d'} \in M_{0}\cup M_3$ for some $d'\le k-1$ such that
 $\bfi_{\calP}(S_0) + \sum_{i\in [d']} \bfi_{\calP}(e_i) \in L$. 

We delete $e_1, \dots, e_{d'}$ from our matching,
thus leaving an unmatched set
$D:=\bigcup_{i\in [d']}e_i\cup S_0$.
Then $\bfi_{\calP}(D)$ is in the integer span of $I$
and $|D| \le k^2(1+\r^{-1})$, so by definition of $C$,
we have $\bfi_{\calP}(D) = \sum_{\bfv\in I} a_{\bfv} \bfv$
as in \eqref{eq:Cmax} with each $|a_{\bfv}|\le C$.
We write each as $a_{\bfv}=b_{\bfv} - c_{\bfv}$ 
such that one of $b_{\bfv}, c_{\bfv}$ 
is a nonnegative integer and the other is zero. 
Thus we have
\begin{equation} \label{eq:replace}
\bfi_{\calP}(D) + \sum_{\bfv\in I} c_{\bfv} \bfv 
= \sum_{\bfv\in I} b_{\bfv} \bfv.
\end{equation}
Let $\F \sub M_2$ contain $c_{\bfv}$ sets
of index $\bfv$ for each $\bfv\in I$.
We delete the family $\F$ of edges from $M_2$,
thus leaving $V(\F)\cup D$ as the unmatched set.

The combinatorial meaning of \eqref{eq:replace} 
is that $V(\F)\cup D$ can be expressed as the
disjoint union over $\bfv\in I$
of $b_{\bfv}$ sets of index $\bfv$.
Thus we regard $V(\F)\cup D$ as the union of 
$|\F|+|D|/k\le kC+k+k/\gamma$
sets $S$ with $\bfi_{\calP}(S)\in I$.
Each such $S$ has at least
$C' - k^2 \ge d' + |\F| + |D|/k$
absorbing sets in $\F_0$,
of which at most $d'$ may be unavailable
due to deleting $e_1, \dots, e_{d'}$,
so we can greedily absorb all such $S$,
thus obtaining a perfect matching of $H$.

For the running time, 
note that we can find $\F_{abs}$ in time $O_k(n^{2^{k+1} k^2})$, 
as $4tk^2=2^{k+1} k^2$.
Then we find $M_2$ and $M_3$ both in time $O_k(n^{k+1})$, 
and $e_1, \dots, e_{d'}$ in time $O_k(n)$.
Finally, it takes constant time to find $b_\bfv$ and $c_\bfv$ 
for all $\bfv\in I$, pick $\F$ and partition $V(\F)\cup D$; 
the absorption can also be done in constant time.
The overall running time is $O_k(n^{2^{k+1} k^2})$.

\section{The extremal case}

In this section we prove Theorem \ref{thm3},
thus completing the proof of Theorem \ref{thm:main}.
We start by recalling the assumptions of the theorem,
which will be assumed throughout the section.
Let $k\ge 3$, suppose $\e>0$ is sufficiently small,
$n>n_0=n_0(k,\e)$ is sufficiently large, 
and let $0 \le c \le \e n/k$ be an integer.
Suppose $H$ is a $k$-graph on $n$ vertices 
with $\delta_{k-1}(H)\ge n/k-c$.
Assume $H\notin \mathcal{H}_{n,k}$ and 
that $H$ is $\e$-extremal, i.e.\ there is 
an independent subset $S\subseteq V(H)$ 
with $|S| \ge (1-\e)\frac{k-1}k n$. 
We will define an algorithm that finds 
either a perfect matching in $H$ 
or a certificate that $H$ has no perfect matching.

\subsection{Proof strategy and notation}

We start by introducing some notation and terminology,
also used throughout the section, which we will motivate
with reference to the strategy of the proof.
We define a partition $(A,B,C)$ of $V(H)$ as follows. 
Let $C$ be any \emph{maximal} independent subset of $V(H)$ that contains $S$. 
Note that one can construct $C$ from $S$ greedily in time $O_k(n^{k})$.
For any $x \in V \sm C$, we write $\deg(x,C)$ for the number 
of edges $e$ of $H$ containing $x$ with $e \sm \{x\} \sub C$.
We let 
\[ 
A = \left\{ x \in V \sm C: 
\deg(x,C) \ge (1-\a) \binom{|C|}{k-1}\right\}, 
\]
where $\alpha=\e^{1/3}$.
The final part of the partition is $B=V(H)\setminus (A\cup C)$.
We will see in Lemma \ref{clm:size} that $B$ is small.

As $C$ is independent, we have $|e \cap C| \le k-1$
for any edge $e$, so we can assume $|C| \le (k-1)n/k$,
as otherwise it cannot be covered by a matching.
Thus $|C| \approx (k-1)n/k$ and $|A| \approx n/k$, 
so the bulk of any perfect matching can be thought 
of as a matching in an auxiliary graph between $A' \sub A$
and a partition of $C' \sub C$ into $(k-1)$-sets.
The main task of the proof will be to efficiently check 
whether we can delete a small matching to reduce to such
an auxiliary graph where a perfect matching can easily 
be found by derandomizing a theorem of Pikhurko
(see Theorem \ref{thm:pik}).

The \emph{slack} of a matching $M$ in $H$ is defined as
\[ s_M := |(A \cup B) \sm V(M)|-(n/k-|M|). \]
As noted above, we can assume
\[ s := s_\es = |A|+|B|-n/k = (k-1)n/k - |C| \ge 0.\]
In general, to use $M$ as part of a perfect matching
we require $s_M \ge 0$. Indeed, writing 
$X = (A \cup B) \sm V(M)$ and $Y = C \sm V(M)$,
if $s_M < 0$ then $|Y| = n - k|M| - |X|
= (k-1) |X| - ks_M > (k-1) |X|$,
so as $Y$ is independent 
there cannot be a perfect matching.

We call a matching $M$ a \emph{cleaner} if $M$ covers $B$
and $s_M$ is non-negative and even. Note that
we impose the final condition to avoid a parity obstruction:
we may only have edges that have odd intersection with $X$,
in which case a perfect matching of $H[X \cup Y]$ requires
$|X|$ and $n/k - |M|$ to have the same parity, namely, $s_M$ is even.
We summarise the above remarks as follows.

\begin{obs} \label{obs:s}
Suppose $M$ is contained in a perfect matching
and $X = (A \cup B) \sm V(M)$. Then:

(i) $s_M \ge 0$, and 
(ii) if all $|e \cap X|$ with $e \in E(H)$
are odd then $s_M$ is even.  
\end{obs}

Next we identify a subset $D$ of vertices in $B$
with relatively high degree within $C$; we will see
in Lemma~\ref{lem:Dext} below that it is easy
to greedily extend any small matching to one covering $D$. 
We list the vertices of $B$ as $v_1,\dots, v_{|B|}$
so that $\deg(v_i,C)$ is a non-increasing sequence.
We find (in linear time) the largest $d > 0$ such that 
$\deg(v_d, C)>(d+2c)(k-1)\binom{|C|}{k-2}$; if no such $d$ exists then let $d=0$.
Let $D:=\{v_1,\dots, v_d\}$ if $d>0$,
or $D:=\emptyset$ otherwise.

We conclude this subsection with some compact notation for
describing the type of a set or edge with respect 
to the partition $(A,B,C)$ of $V(H)$.
We say that $S$ is an $A^i B^j C^l$ set
if $|S\cap A|=i$, $|S\cap B|=j$ and $|S\cap C|=l$;
we also say $S$ has the form $A^i B^j C^l$.
If any index is $0$, we omit it, e.g.\ we write 
$B^jC^{l}$ instead of $A^0B^jC^{l}$.
If $S$ is an edge of $H$ we call it an $A^i B^j C^l$ edge.
We also call it an $(i+j)$-edge. 
We write $N_H(v, A^i B^j C^l)$ for the $(k-1)$-sets in $N_H(v)$ 
of form $A^i B^j C^l$, namely,
\[
N_H(v, A^i B^j C^l) := \{ S\in N_H(v): |S\cap A|=i, |S\cap B|=j, |S\cap C|=l\}
\]
and also
$\deg(v,A^i B^j C^l) = |N_H(v, A^i B^j C^l)|$.

\subsection{The algorithm and proof modulo claims}

Now we define our algorithm 
(it refers to some claims stated below).
The notation is as in the previous subsection; we also let
\[ t = n/k - |A| = |B| - s
\ \text{ and } \ d=|D|. \]
 
\begin{procedure}[h]
  \caption{PerfectMatchingEXT()}
  \label{alg:determine}
  \SetKwInOut{Input}{Data}
  \SetKwInOut{Output}{Output}
  \Input{An $\e$-extremal $n$-vertex $k$-graph $H\notin \mathcal H_{n,k}$
  with $\delta_{k-1}(H)\ge n/k-c$, 
  an independent set $S$ with $|S|\ge (1-\e)\frac{k-1}{k}n$,
  and sets $A$, $B$, $C$, $D$ as defined above.}
  \Output{A perfect matching in $H$ or a certificate that none exists.} 

  \If {$n<n_0$} {
      Examine every set of $n/k$ edges in $H$, and halt with appropriate
output.   }

  \If {$|C|>(k-1)n/k$} {
      Output ``no PM'' and $C$, and halt.  }

  \If {$s=|A|+|B|-n/k$ is odd and 
  $H$ contains no $j$-edge with $j\le s+1$ even} 
  {Output ``no PM'' and $H$, and halt.}
  
  \eIf {there is no matching $M_0$ of size $\max\{t-d,0\}$ 
  in $H[(B\setminus D)\cup C]$} 
  {Output ``no PM'' and $(B\setminus D)\cup C$, and halt.}
  {Output ``PM'' \;}
     
 Use one of Lemmas \ref{lem:last}, \ref{lem:main-}
 or \ref{lem:main} to find a cleaner $M$,
 then Lemma \ref{lem:finish} to find a perfect matching.
\end{procedure} 

Next we state five claims needed for the proof
of Theorem \ref{thm3} (the proofs are deferred to
later in the section). The first claim is required
so that the search for $M_0$ is feasible, 
the second shows that it is sufficient to find a cleaner, 
and the others show how to find a cleaner in various cases.

\begin{claim}\label{lem:t-d}
We have $t-d \le 2c$, 
and $M_0$ can be found in time $O(n^k +2^{O(kc)})$.
\end{claim}

\begin{claim}\label{lem:finish}
If $H$ has a cleaner $M$ then $H \sm V(M)$ has a perfect matching, 
which can be found in time $O(n^{4(k-1)})$.
\end{claim}

\begin{claim}\label{lem:last}
Suppose $H$ contains no $j$-edge for all even $0\le j\le k$ 
and $H\notin \mathcal{H}_{n,k}$.
If Procedure~\ref{alg:determine} outputs ``PM'' 
then we can find a cleaner in time $O_k(n^k)$.
\end{claim}

\begin{claim}\label{lem:main-}
Suppose $s=0$ or that
$H$ has an $i$-edge for some even $i \in [2, k]$
but does not have any $ABC^{k-2}$ edge.
If Procedure~\ref{alg:determine} outputs ``PM'' 
then we can find a cleaner in time $O_k(n^k)$.
\end{claim}

\begin{claim}\label{lem:main}
Suppose $s>0$ and $H$ has some $ABC^{k-2}$ edge $e_0$.
If Procedure~\ref{alg:determine} outputs ``PM'' 
then we can find a cleaner in time $O_k(n^k)$.
\end{claim}

We remark that when Procedure~\ref{alg:determine} outputs ``no PM'' (for lines 2, 4, 6 or 8),  it also outputs a certificate.
For example, when it outputs the set $C$ (line 4), the set $C$ is an independent set of $H$ of size larger than $(k-1)n/k$, implying that $H$ has no perfect matching.
When it outputs $H$ (line 6), we have that $s=|A|+|B|-n/k$ is odd and 
  $H$ contains no $j$-edge with $j\le s+1$ even.
Since $s=|A|+|B|-n/k$ is odd, any perfect matching of $H$ must contain an edge that contains an even number of vertices in $A\cup B$.
However, such an edge $E$ must contain at least $s+3$ vertices from $A\cup B$, meaning that the slack of the matching $\{E\}$ is at most $(|A|+|B|-(s+3))-(n/k-1) = -2 <0$, and thus $E$ is not contained in any perfect matching, a contradiction.
Finally, for line 8, note that by the definition of $t$, if $H$ has a perfect matching, $H[B\cup C]$ must contain a matching of size $n/k-|A|=t$, and thus $H[(B\setminus D)\cup C]$ must contain a matching of size $\max\{t-|D|, 0\}$.

We conclude this subsection by assuming these claims
and proving Theorem \ref{thm3}. First we make the
following observation which will henceforth
be used without comment.

\begin{obs}
If a matching $M$ has $n_i$ $i$-edges then 
$s_M = s - \sum_{i \ge 1} n_i (i-1)$. 
\end{obs}

\begin{proof}[Proof of Theorem \ref{thm3}]
We need to show that Procedure~\ref{alg:determine} either finds a perfect matching in $H$ 
or a certificate that none exists, in time $O(n^{4(k-1)})+O_k(n^k+ 2^{O(kc)})=O(n^{4k-4}+ 2^{O(kc)})$.

For the running time, we find $C$ as an arbitrary 
maximal independent set containing $S$ in time $O(n^{k})$, 
and then find $A$ and $B$ by determining their degrees to $C$ 
in time $O_k(n^{k})$. Next we sort the degrees of vertices from $B$ 
in time $O(n\log n)$ and find $d$ and $D$ in linear time.
The first two tests runs in constant time 
and the third in time $O_k(n^k)$. By Claim \ref{lem:t-d}, 
the search for $M_0$ takes time $O(n^k+2^{O(kc)})$.
If the output is ``PM'' then the conclusion of the algorithm
finds a perfect matching in time $O(n^{4(k-1)})$
(assuming the claims).

It remains to show that there is no perfect matching
if the algorithm outputs ``no PM''. The correctness of the
first three tests follows from Observation \ref{obs:s}.
Indeed, for the third test, let $M=\{e\}$ for any $j$-edge $e$
with $j\ge s+3$ even.
Then $s_M = s - (j-1) \le -2$, so no matching containing $e$
is perfect.
On the other hand, as $s$ is odd, any matching consisting of odd edges only 
cannot be perfect either.
For the final test, note that in any matching $A\cup D$ 
can be incident to at most $|A \cup D| = n/k-t+d$ edges,
so any perfect matching must contain a matching of size 
$\max\{t-d,0\}$ completely within $(B\setminus D)\cup C$.
\end{proof}

\subsection{Using a cleaner}

In this subsection we prove Claim \ref{lem:finish},
which shows how to find a perfect matching 
assuming that there is a cleaner.
We start by estimating the sizes of the parts $A$, $B$, $C$;
as discussed above, $B$ is small, 
$|A| \approx n/k$ and $|C| \approx n - n/k$.

\begin{claim}\label{clm:size}
$|A|\ge n/k-\a^2 n$, $|B|\le \a^2 n$ 
and $(1-{\e})\frac {(k-1)n}k\le |C|\le \frac {(k-1)n}k$.
\end{claim}

\begin{proof}
The upper bound on $C$ follows from Observation \ref{obs:s}
and the lower bound from our assumptions for Theorem \ref{thm3},
which give $|C| \ge |S| \ge (1-{\e})(k-1)n/k$.
As $|A|+|B|+|C|=n$, we have
\begin{equation}\label{eq:ab}
0 \le s = |A|+|B|-n/k =n-n/k-|C| \le \e (k-1)n/k.
\end{equation}
By the definitions of $A$ and $B$, we have
\[
\left(\frac{n}k -c \right) \binom{|C|}{k-1} 
\le \sum_{x \in A \cup B} \deg(x,C) \le (1-\a)\binom{|C|}{k-1} |B| + \binom{|C|}{k-1} |A|.
\]
We deduce $n/k -c \le |A|+|B|-\a |B|$, 
so $\a |B| \le |A|+|B| - n/k+c\le \e n$ by \eqref{eq:ab} 
and $c \le \e n/k$.
Thus $|B|\le \a^2 n$ and by \eqref{eq:ab} again, 
$|A|\ge n/k-|B|\ge n/k-\a^2 n$. 
\end{proof}

We also require the following algorithmic version 
of a special case of a result of Pikhurko \cite{Pik}, 
concerning perfect matchings 
in a $k$-graph $H$ that is $k$-partite, 
i.e.\ $V(H)$ has a partition $(V_1,\dots,V_k)$ 
so that every edge intersects all $k$ parts.
For $S \sub [k] = \{1,\dots,k\}$ we write
$\delta_S(H)$ for the minimum degree 
$\deg_H(\{v_i: i \in S\})$ of any set
consisting of one vertex in each 
of the parts $(V_i: i \in S)$.

\begin{theorem}\cite[Theorem 3]{Pik}\label{thm:pik}
Suppose $1/n\ll \gamma \ll 1/k$.
Let $H$ be a $k$-partite $k$-graph with
parts $V_1,\dots,V_k$ each of size $n$.
Suppose $\delta_{\{1\}}(H) \ge (1-\gamma) n^{k-1}$ 
and $\delta_{[k]\setminus\{1\}}(H)\ge (1-\gamma) n$.
Then there is an algorithm that finds a perfect matching in $H$ 
in time $O(n^{4(k-1)})$.
\end{theorem}

As Pikhurko's proof is probabilistic, we give 
an alternative derandomised proof via Claim~\ref{lem:deran}
(see also \cite{HZ3} for a similar proof).

\begin{proof}
We will apply Lemma~\ref{lem:deran} to the auxiliary graph $G$
with parts $U=V_1$ and $W=\prod_{i=2}^k V_i$, 
whose edges consist of all $\{u,S\}$
with $u \in U$, $S \in W$ and $\{u\} \cup S \in E(H)$,
and all $\{S,S'\}$ with $S,S' \in W$ and $S \cap S' \ne \es$.
According to the notation of Lemma~\ref{lem:deran}
we have $M=n$, $N=n^{k-1}$ and
$m = |G[W]| \le (k-1)n^{2k-3}$. 
For any $u \in U$ we have
$\deg_G(u, W) \ge \delta_{\{1\}}(H) \ge (1-\gamma)n^{k-1}$,
so we can take $\beta = 1 - \gamma$.
To satisfy the conditions of Lemma~\ref{lem:deran}
we let $\tau := 1/4$ and $r:=4\gamma n$.
Then for large $n$ we have\[
\exp \left(\frac{\tau^2 r}{3(1-\gamma)}\right) 
\ge \exp \left( \frac{\gamma n}{12} \right) \ge 8n = 8M.
\]
By Lemma~\ref{lem:deran}, in time $O(N^4+MN^3) = O(n^{4(k-1)})$, 
we can find an independent set $R$ in $G[W]$ 
such that $(1-\nu)r\le |R|\le r$ 
and $|N_G(u)\cap R| \ge(\beta - \tau - \nu)r$ for all $u\in U$, 
where \[
\nu = \frac{2m r}{N^2}\le \frac{8(k-1)\gamma n^{2k-2}}{n^{2k-2}}
= {8(k-1)\gamma}.
\] 
Thus $R$ is a collection of disjoint $(k-1)$-sets
of the form $v_2 \dots v_k$ with each $v_i \in V_i$
with $(1-\nu)4\gamma n \le |R| \le 4\gamma n$
and $|N_G(u)\cap R| \ge(1-8k\gamma - 1/4)r>r/2=2\gamma n$ 
for all $u\in V_1$, for small $\r$.

Next we partition $(\bigcup_{2\le i\le k}V_i)\setminus V(R)$ 
arbitrarily into $n-|R|$ disjoint $(k-1)$-sets 
$S_1, S_2, \dots, S_{n-|R|}$ each
having one vertex in each of $V_2,\dots,V_k$.
We consider the bipartite subgraph $G'$ of $G$ 
induced by $V_1$ and $\{S_1, \dots, S_{n-|R|}\}$.

We claim that the set $V_1'$ of vertices $v \in V_1$
with $\deg_{G'}(v) < (n-|R|)/2$ has size at most $2\gamma n$.
To see this, first note that $|E(G')|\ge (1-\gamma)n(n-|R|)$
as $\delta_{[k]\setminus\{1\}}(H) \ge (1-\gamma) n$.
Thus $G'$ has at most $\gamma n(n-|R|)$ non-edges,
whereas every vertex in $V_1'$ is in at least 
$(n-|R|)/2$ non-edges of $G'$, so the claim holds.

We now proceed as follows.
\begin{enumerate}
\item Greedily match the vertices of $V_1'$ with $(k-1)$-sets in $R$. 
This is possible because each vertex $v$ forms an edge 
with more than $2\gamma n$ sets in $R$, and $|V_1'|\le 2\gamma n$.
\item Greedily match the remaining sets of $R$ 
with vertices in $V_1\setminus V_1'$. 
This is possible as $|R| \le 4\gamma n$
and each set in $R$ forms an edge 
with at least $(1-\gamma)n$ vertices in $V_1$.
\item Let $V_1''$ be the set of uncovered vertices in $V_1$
and $G''$ be the bipartite subgraph of $G$
induced by $V_1''$ and $\{S_1, \dots, S_{n-|R|}\}$.
We note that both parts of $G''$ have size $n-|R|$
and $G''$ has minimum degree at least $(n-|R|)/2$.
Thus $G''$ has a perfect matching by Hall's theorem,
which can be found in time $O(n^3)$ by the Hungarian algorithm.
\end{enumerate}
Clearly the union of the matchings above gives a perfect matching in $H$.
\end{proof}

\begin{proof}[Proof of Claim~\ref{lem:finish}]
Suppose $H$ has a cleaner, i.e.\ a matching $M$
that covers $B$ such that
$s_M = |(A \cup B) \sm V(M)|-(n/k-|M|)$
is non-negative and even. 
We will show that $H \sm V(M)$ has a perfect matching, 
which can be found in time $O(n^{4(k-1)})$.

We start by finding a matching $M_*$
of size $2\alpha^2 n$ with $V(M_*) \cap V(M) = \es$ 
consisting only of $A^2 C^{k-2}$ edges
or only of $A^3 C^{k-3}$ edges. To do so,
we apply a greedy algorithm to choose edges one by one,
which we call `good' if they avoid $V(M)$ and all previous choices, 
until we obtain $2\alpha^2 n$ edges of one of the required forms.
At each step, we consider any good set $S$ consisting
of $2$ vertices in $A$ and $k-3$ vertices in $C$.
As $\deg_H(S) \ge n/k - c > |B| + k|M| + 2k\alpha^2 n$
we can complete $S$ to a good edge of the form 
$A^2C^{k-2}$ or $A^3C^{k-3}$. Thus the algorithm to find $M_*$
can be completed, and clearly takes time $O(n^2)$.

We obtain $M'$ from $M$ by adding some edges of $M_*$:
we add $s_M$ edges if they have the form $A^2C^{k-2}$ 
or $s_M/2$ edges if they have the form $A^3C^{k-3}$.
In both cases we obtain $s_{M'} = 0$,
i.e.\ $|C'| = (k-1)|A'|$ where 
$A' = A \sm V(M')$ and $C' = C \sm V(M')$.
We partition $C'$ arbitrarily into $k-1$ parts
$C^1, C^2,\dots, C^{k-1}$ each of size $m:=|A'|$.
We note by \eqref{eq:ab} that $m \ge |A| - k|M'| 
\ge |A|-k(|B|+s_M) > n/k - 2k\a^2 n$.

To find a perfect matching, it suffices to show
that Theorem~\ref{thm:pik} applies to the $k$-partite 
sub-$k$-graph $H'$ of $H$ with parts
$A', C^1, C^2,\dots, C^{k-1}$.

To bound $\delta_{[k]\setminus\{1\}}(H')$
we consider any set $S$ formed by $k-1$ vertices
$v_i \in C^i$ for $i \in [k-1]$. As $C$ is independent, 
the number of non-neighbours of $S$ in $A \cup B$
is at most
\[
|A| + |B| - {n}/{k}+c \le \e \frac {(k-1)n}k +c\le k\e m,
\]
where we use \eqref{eq:ab} in the first inequality and the last inequality follows from $m=|A'|\ge n/k - 2k\a^2 n>\frac{k-1}{k^2}n$.
We deduce $\delta_{[k]\setminus\{1\}}(H')\ge m - k\e m= (1-k\e)m$.

To bound $\delta_{\{1\}}(H')$, we consider any $v \in A' \sub A$,
and note that by definition of $A$, the number of $(k-1)$-sets 
$S \sub C$ with $S \notin N_H(v)$ is at most
\[
\a \frac{|C|^{k-1}}{(k-1)!}\le \a \frac{\left(\frac{k-1}{k}n \right)^{k-1}}{(k-1)!}\le \a \frac{(km)^{k-1}}{(k-1)!}= \a c_k m^{k-1},
\]
where $c_k=\frac{k^{k-1}}{(k-1)!}$.
This implies $\delta_{\{1\}}(H') \ge (1-\a c_k) m^{k-1}$. 

By Theorem \ref{thm:pik} with $\gamma=\a c_k$, 
we find a perfect matching in $H'$ in time $O(n^{4(k-1)})$.
\end{proof}

\subsection{Greedy extension}

In this subsection we gather various results 
that will be used in the following subsection
for greedy extension of matchings
during the construction of a cleaner.

\begin{claim}\label{lem:M2}
Let $B'\subseteq B$ and $X \subseteq A \cup C$ 
with $|X| \le k|B|$.
Then in time $O(n^2)$ we can find a matching $M_2$ 
with $V(M_2)\cap X=\emptyset$
that covers all vertices of $B'$ by edges 
of the form $ABC^{k-2}$ or $BC^{k-1}$.
\end{claim}

\begin{proof}
For each $v\in B'$ in turn
we pick $k-2$ arbitrary vertices 
from $C\setminus X$,
and an uncovered vertex in $V\setminus (B\cup X)$ 
to complete an edge, which we add to $M_2$. 
Since $|B|+|X|\le 2k|B|<n/k-c\le \delta_{k-1}(H)$, 
such an edge always exists.
\end{proof}

Next we recall that $D$ consists of $d$ vertices
each having degree in $C$ larger than 
$(d+c)(k-1)\binom{|C|}{k-2}$.
To see the following claim,
we consider a greedy algorithm for constructing $M$,
and note that at each step at most $(k-1)(d-1+c+1)$ vertices in $C$ 
are already chosen, so the number of unavailable $(k-1)$-sets in $C$
is at most $(d+c)(k-1)\binom{|C|}{k-2}$.
 
\begin{claim}\label{lem:Dext}
Given any set $C_*\subseteq C$ of size at most $(c+1)(k-1)$, we can find a matching of $M$ of size $d$ such that $D\subseteq V(M)\subseteq D\cup (C\setminus C_*)$ in time $O(n^{k})$.
\end{claim}

Now we give a degree bound that will be
used to greedily cover $B \sm D$.

\begin{claim}\label{lem:D}
If $D\subsetneq B$ then $\deg(v, AC^{k-2})>k|B||A|\binom{|C|}{k-3}$ 
for every vertex $v\in B\setminus D$.
\end{claim}

\begin{proof}
Suppose for contradiction that there is $v\in B\setminus D$
with $\deg(v, AC^{k-2}) \le k|B||A|\binom{|C|}{k-3}$. 
By considering the degree of $v$ 
together with each $(k-2)$-set in $C$ we have
\[
\deg(v, AC^{k-2}) + \deg(v, BC^{k-2}) + (k-1)\deg(v, C^{k-1}) 
\ge \binom{|C|}{k-2} \left(\frac nk -c \right).
\]
As $|B|\le \alpha n^2$ by Claim~\ref{clm:size} 
and $\alpha$ is small we deduce
\[
(k-1)\deg(v, C^{k-1}) 
\ge \binom{|C|}{k-2} \left(\frac nk -c \right) - k|B||A|\binom{|C|}{k-3} - |B|\binom{|C|}{k-2} 
\ge k^2(|B|+2c)\binom{|C|}{k-2}.
\]
However, this contradicts the definition of $D$, so the claim holds.
\end{proof}

We conclude this subsection by establishing 
the lower bound $d=|D| \ge t-c$
and describing the algorithm to find $M_0$.
 
\begin{proof}[Proof of Claim \ref{lem:t-d}]
Note that for any vertex $b\in B\setminus D$, if $d>0$ then $\deg(b, C)\le (d+2c)(k-1)\binom{|C|}{k-2}$; if $d=0$ then $\deg(b, C)\le (1+2c)(k-1)\binom{|C|}{k-2}$.
So it always holds that $\deg(b, C)\le (d+1+2c)(k-1)\binom{|C|}{k-2}$.
Recall that  
(i) $\delta_{k-1}(H[B\cup C])\ge t-c$ and $C$ is independent,
 (ii) $|B|\le \alpha^2 n$ and $|C|\le (k-1)n/k$, and (iii) $\deg(b, c)\le (1-\alpha)\binom{|C|}{k-1}$ for every $b\in B$.
Putting these together, we get
\[
\begin{aligned}
\binom{|C|}{k-1}(t-c) &\le \sum_{b\in B} \deg(b, C) \\
&\le d(1-\alpha)\binom{|C|}{k-1} + |B| (d+1+2c)(k-1)\binom{|C|}{k-2} 
\le (d(1-\alpha)+\alpha(d+2c+1)) \binom{|C|}{k-1},
\end{aligned}
\]
that is, $t-c\le d+2\alpha c +\alpha$.
This implies that $t-d\le 2c$ as $\alpha<1/3$.

Now assume $t-d>0$, and we need to find a matching $M_0$ of size $t-d\le 2c$ in $H[(B\setminus D)\cup C]$, if one exists.
This is done by Theorem~\ref{thm:Fellows} in time $O(n^k+2^{O(kc)})$.
\end{proof}

\subsection{Finding a cleaner}

We conclude this section by proving the three claims
that find a cleaner according to the various cases
of the proof of Theorem \ref{thm3}.

\begin{proof}[Proof of Claim \ref{lem:last}]
Assume that $H$ contains no even edges 
and $H\notin \mathcal{H}_{n,k}$.
Then the slack $s=|A|+|B|-n/k$ must be even. 
We apply Claim~\ref{lem:M2} to find a matching $M$ covering $B$,
which must consist of $BC^{k-1}$ edges, as there is no 2-edge.
This preserves the slack, i.e.\ $s'=s=0$,
so $M$ is a cleaner, as required.
\end{proof}

\begin{proof}[Proof of Claim \ref{lem:main-}]
We start with the case $s=0$.
Procedure~\ref{alg:determine} finds a matching $M_0$ 
of size $t-d=|B|-|D|$ in $(B\sm D)\cup C$,
as otherwise it would output ``no PM''.
By Claim~\ref{lem:Dext}, we can enlarge $M_0$ to a matching $M$ 
of size $|B|$ covering $B$ and $(k-1)|B|$ vertices in $C$.
Thus we preserve the slack, i.e.\ $s'=s=0$,
so $M$ is a cleaner, as required.

It remains to consider $s>0$.
By assumption, $H$ has an $i$-edge 
for some even $i \in [2, k]$
and does not have any $ABC^{k-2}$ edge.
Thus whenever we apply Claim~\ref{lem:M2}, 
we always obtain edges of the form $BC^{k-1}$.
If $s$ is even then we can simply apply Claim~\ref{lem:M2} 
to construct a matching $M$ covering $B$ 
with edges of the form $BC^{k-1}$.
Indeed, this preserves the slack, i.e.\ $s'=s$,
so $M$ is a cleaner, as required.

Finally, we can assume that $s$ is odd.
By assumption, there is an $i$-edge $e_0$ with $i$ even.
We fix any $e_0$ that minimises $i$.
Note that $s \ge i-1$, otherwise
Procedure~\ref{alg:determine} would output ``no PM''.
We apply Claim~\ref{lem:M2} 
to construct a matching $M$ consisting of $e_0$
and $|B \sm e_0|$ edges covering $B \sm e_0$ 
with edges of the form $BC^{k-1}$.
Thus we obtain slack
$s'=|A \sm V(M)|-(n/k-|M|)=s-i+1$,
which is even and non-negative,
so again $M$ is a cleaner, as required.
\end{proof}

\begin{proof}[Proof of Claim \ref{lem:main}]
Suppose $s>0$ and
$H$ has some $ABC^{k-2}$ edge $e_0$,
i.e.\ $|e_0 \cap A|=|e_0 \cap B|=1$,
and Procedure~\ref{alg:determine} outputs ``PM''.
We write $e_0 \cap B = \{x\}$.
Our proof will use two different strategies
for finding a cleaner that in combination
cover all possible cases. The main strategy 
(which applies to all but one case)
finds in time $O_k(n^k)$ a matching 
$M = M_1 \cup M_2$ covering $B$, such that, writing
$s_2 := s_{M_2} = |(A \cup B) \sm V(M_2)| - (n/k-|M_2|)$,
\begin{enumerate}
\item every edge $e$ of $M$ has 
$|e \cap B| = 1$ and $|e \cap A| \le 1$, 
\item $|M_1| \ge t+1 = n/k - |A| + 1 = |B|-s+1$,
\item the number of edges of $M_1$ that intersect $A$
is $(s_2 \text{ mod } 2) \in \{0,1\}$.
\end{enumerate}
Note that if we find such a matching $M$ 
then $|M_2| = |M|-|M_1| \le |B|-(t+1) = s-1$ 
and $s_2 = s + |M_2| - |(A \cup B) \cap V(M_2)|
\ge s - |M_2| \ge 1$, 
so $s_M = s_2 - (s_2 \text{ mod } 2)$
is even and non-negative,
i.e.\ $M$ is a cleaner, as required.
We will apply the main strategy to $3$ cases below.

\medskip

\nib{Case 1.} Suppose $D=B$. By Claim~\ref{lem:Dext},
we can find a matching $M'$ covering $B$ by $BC^{k-1}$ edges
such that $V(M')$ is disjoint from $e_0 \sm B$.
Then $|M'|=|B|=t+s \ge t+1$, as $s>0$.
Let $e_x$ be the edge of $M'$ containing $x$
(the vertex in $e_0 \cap B$).
We let $M_2 = \es$, so $s_2=s=|A|+|B|-n/k$,
and let $M_1=M'$ if $s_2$ is even
or otherwise let $M_1 = (M' \sm \{e_x\}) \cup \{e_0\}$.
Then $M$ satisfies (1-3).

\medskip

%

\nib{Case 2.} Suppose $t \le |D| < |B|$.
Fix any $v \in B \sm D$. By Claim~\ref{lem:D}, we can find
an $ABC^{k-2}$ edge $e'$ containing $v$.
Moreover, as $C$ is a maximal independent set, we can find 
a $BC^{k-1}$ edge $e_v$ containing $v$.
Next, by Claim~\ref{lem:Dext} we can find a matching $M'$
covering $D$ by $BC^{k-1}$ edges with 
$V(M') \cap (e_v \cup e') = \es$.
Then by Claim~\ref{lem:M2} we find a matching $M_2$
covering $B \sm (D \cup \{v\})$ with 
$V(M_2) \cap (e_v \cup e' \cup V(M')) = \es$.
We let $M_1 = M' \cup \{e_v\}$ if $s_2$ is even
or $M_1 = M' \cup \{e'\}$ otherwise.
Then $M$ satisfies (1-3).

\medskip

\nib{Case 3.} Suppose $|D|=d<t$ and there is 
a $BC^{k-1}$ edge $e_*$ disjoint from $D \cup V(M_0)$,
where $M_0$ is the matching of size $t-d$ in $(B\sm D)\cup C$
found by Procedure~\ref{alg:determine}.
We write $e_* \cap B = \{v\}$.
By Claim~\ref{lem:Dext}, we can find $BC^{k-1}$ edges 
covering $D$ that extend $M_0 \cup \{e_*\}$ 
to a matching $M'$ of size $t+1$.
Since $v\notin D$, by Claim~\ref{lem:D} we have 
$\deg(v, AC^{k-2})>k|B||A|\binom{|C|}{k-3}
\ge k|M'||A|\binom{|C|}{k-3}$, so we can find 
an $ABC^{k-2}$ edge $e_v$ containing $v$ 
such that $e_v \cap V(M')=\{v\}$.
By Claim~\ref{lem:M2} we can find $M_2$
covering $B \sm V(M')$ such that
$V(M_2)\cap (V(M')\cup e_v)=\es$.
We let $M_1=M'$ if $s_2$ is even
or $M_1 = (M' \sm \{e_*\}) \cup \{e_v\}$ otherwise.  
Then $M$ satisfies (1-3).

\medskip

Finally, we describe the second strategy for finding a cleaner,
which will complete the proof when none of the above 3 cases apply,
i.e.\ when $|D|=d<t$ and there is no
$BC^{k-1}$ edge $e_*$ disjoint from $D \cup V(M_0)$.
By Claim~\ref{lem:Dext}, we can find $BC^{k-1}$ edges 
covering $D$ that extend $M_0$ to a matching $M'$ of size $t$.

Then we claim that we can find $ABC^{k-2}$ edges
covering $B \sm V(M')$ that extend $M'$ to a matching $M$.
To see this, we apply a greedy algorithm, where in each step, 
to cover some vertex $x$ of $B \sm V(M')$ we consider any 
$BC^{k-2}$ set $S$ containing $x$ disjoint from all previous edges.
By our assumptions for this case, $N(S)$ is disjoint from $C$,
so as $|N(S)| \ge n/k - c > k|B|$ we can complete $S$
to an $ABC^{k-2}$ edge as required, 
so we can construct $M$ as claimed.

Writing $t' = |V(M_0) \cap B|$, we note that 
$M_0$ reduces the slack by $t'-|M_0|=t'-t+d$,
and $M' \sm M_0$ keeps it unchanged.
Finally, $M\sm M'$ reduces the slack by $|M \sm M'| = |B|-d-t'$.
Therefore $s_M = s - (t'-t+d) - (|B|-t'-d) = 0$, 
so $M$ is a cleaner.
\end{proof}

\section{Concluding remarks}

In this paper, we showed that $\DPM_k(n,m)$ 
is in P for $m\ge n/k-c$ for $c=O_k(\log n)$.
For simplicity  of analysis we did not attempt 
to optimise the exponent in the running time.
We remark that some improvements can be obtained
by merging `transferrals' and working 
with two partitions as in~\cite{Han14_poly},
and/or settling for just testing for a perfect matching
rather than actually finding one.

{
For matchings that are not necessarily perfect (Theorem~\ref{coro}), note that our proof of Theorem~\ref{coro} indeed shows that $H$ with $\delta_{k-1}(H')\ge m-c$ contains a matching of size $m$ if and only if the auxiliary graph $H'$ with $\delta_{k-1}(H')\ge |H'|/k-c$ contains a perfect matching.
Moreover, when $c = \Omega(n)$, it is shown in~\cite{Szy13} that the NP-Complete problem $\DPM_k(n,0)$ admits a polynomial reduction to $\DPM_k(n,n/k-c)$.
Therefore, under the Exponential-Time Hypothesis, there is no algorithm that solves the \textsc{$k$-Graph Matching Above Minimum Codegree} problem for $c = \Omega(n)$ within time $2^{o(c)}n^{g(k)}$.
}

As observed in~\cite{Han14_poly}, the argument in \cite{Szy13} 
actually shows that $\DPM_k(n, n/k - n^\e)$ 
is NP-complete for any $\e>0$.
Thus the complexity remains unknown 
for $n/k - n^{\e}\le m< n/k-\omega(\log n)$.
Our algorithm is valid for any $c=o(n)$,
but our proof only gives polynomial running time when $c=O_k(\log n)$.
The bottleneck to improving this comes from the extremal case, 
where we showed that the existence of a perfect matching in $H$ 
is equivalent to that of a matching of size $t:=n/k-|A|$ in $H[B\cup C]$.
Thus we have a reduction to the following potentially simpler problem,
which we believe has independent interest.

\begin{problem}\label{prob}
Suppose $1/n\ll \e \ll 1/k$.
Let $X$, $Y$ be disjoint sets with $|X|= \e n$ and $|Y|=n$.
Let $c=c(n)$ and $t=t(n)$ with $c\le t\le |X|$.
Let $H$ be a $k$-graph on $X \cup Y$ 
with $\delta_{k-1}(H)\ge t-c$ 
such that $H[Y]$ is independent. 
What is the complexity of deciding
the existence of a matching of size $t$ in $H$?
\end{problem}

\section{Acknowledgement}
We thank the anonymous referees for detailed comments that have substantially improved the presentation of this paper.

\bibliography{Bibref}

\end{document}